\documentclass{article}
\usepackage[dvips]{graphics}
\usepackage{graphicx}
\usepackage{amscd}
\usepackage{amsmath}
\usepackage{amssymb}
\usepackage[mathscr]{euscript}
\usepackage[all]{xy}

\usepackage{comment}
\numberwithin{equation}{section}
\newtheorem{thm}{Theorem}[section]
\newtheorem{theorem}[thm]{Theorem}
\newtheorem{cor}[thm]{Corollary}
\newtheorem{prop}[thm]{Proposition}
\newtheorem{lem}[thm]{Lemma}

\newcommand{\alp}{\alpha}

\newcommand{\ep}{\epsilon}
\newcommand{\lam}{\lambda}
\newcommand{\Lam}{\varLambda}
\newcommand{\Ome}{\varOmega}
\newcommand{\bll}{\bigl|}
\newcommand{\brl}{\bigr|}
\newcommand{\FS}{\mathrm{FS}}
\newcommand{\gam}{\gamma}
\newcommand{\Gam}{\varGamma}
\newcommand{\II}{{\cal I}}
\newcommand{\inte}{{\textstyle \int}}

\newcommand{\Sym}{{{\mathfrak S}}}
\newcommand{\Sn}{{{\mathfrak S}_n}}

\newcommand{\Ind}{{\mathrm{Ind}}}

\newcommand{\PP}{ {\mathcal P}}

\newcommand{\Res}{{\mathrm{Res}}}
\newcommand{\STab}{{\mathrm{STab}}}
\newcommand{\tb}{{\overline{\otimes}}}
\newcommand{\tr}{{\mathrm{Tr}}}
\newcommand{\id}{{\mathrm{id}}}
\newcommand{\cu}{\varepsilon}

\newcommand{\CCC}{{\cal C}}

\newcommand{\FFF}{{\cal F}}
\newcommand{\C}{{\Bbb C}}

\newcommand{\Z}{{\Bbb Z}}

\newcommand{\FGX}{{\cal F} (G,X)}

\newcommand{\uob}{{\boldsymbol{1}}}
\newcommand{\tra}{{\sf T}}

\newcommand{\KXG}{\C^X\! \rtimes  G}

\newcommand{\hako}[4]{\begin{xy}
(0,2.5)="A", (5,2.5)="B",
(0,-2.5)="C", (5,-2.5)="D", 
\ar@{-}"A";"B"^{#1}
\ar@{-}"C";"D"_{#2} 
\ar@{-}"A";"C"_{#3} 
\ar@{-}"B";"D"^{#4}
\end{xy}}
\newcommand{\hakoko}[6]{\begin{xy}
(0,2.5)="A", (5,2.5)="B",
(0,-2.5)="C", (5,-2.5)="D", 
\ar@{-}"A";"B"^{#1}
\ar@{-}"C";"D"_{#2} 
\ar@{-}"A";"C"_{\genfrac{}{}{0pt}{1}{#3}{#4}} 
\ar@{-}"B";"D"^{\genfrac{}{}{0pt}{1}{#5}{#6}} 
\end{xy}}
\newcommand{\bou}[3]{\begin{xy}
(0,2.5)="A", (5,2.5)="B",
(0,-2.5)="C", (5,-2.5)="D", 
(2.5,2.5)="E", (2.5,-2.5)="F",
\ar@{}"A";"B"^{#1}
\ar@{}"C";"D"_{#2} 
\ar@{-}"E";"F"^{#3}
\end{xy}}

\begin{document}
%
%
%\newcommand{\kasen}[1]{\underline{\text{#1}}}
%\newcommand{\onto}{\twoheadrightarrow}

%

%
%				
%\title{Bialgebroids, % non-commutative 
%A-B torsors and equivalences of equivariant vector bundles %(preliminary  version)}

%\title{Bialgebroids, % non-commutative 
%A-B torsors and equivalences of equivariant vector bundles}
\title{Weak Hopf algebras and the distribution of involutions in symmetric groups}

\author{Takahiro Hayashi \quad}
%\email{hayashi@@math.nagoya-u.ac.jp}
\date{$\qquad$\\
%\date{%$\qquad$ 
Graduate School of Mathematics, 
Nagoya University, \\ %$\quad\;$
Chikusa-ku, Nagoya 464-8602, Japan\\
E-mail: hayashi@math.nagoya-u.ac.jp}
\maketitle

\begin{abstract}
By computing Frobenius-Schur indicators of modules of certain weak Hopf algebras, 
we give a formula for the number of involutions in symmetric groups, which are contained in a given coset with respect to a given Young subgroup.   
\end{abstract}

%\newpage
\section{Introduction}

Let $\Sn$ be the symmetric group of $n$-letters.
Then we have the following classical identity in combinatorial representation theory:
\begin{equation}
\label{InvSn}
\bll \{ a \in \Sym_n \,|\, a^2 = 1 \}\brl
=
 \bll\STab (n)\brl,
\end{equation}
where $\STab (n)$ denotes the set of all standard tableaux of size $n$. 
Besides a proof based on RSK correspondence, there is a proof of this identity,
which is based on Frobenius-Schur indicators. 
Let $G$ be an arbitrary finite group. Then the {\it $r$-th root number function} 
$R_{G}^r (a) := \bll \{ c \in G\,|\, c^r = a \}\bll$
is given by
\begin{equation}
\label{RGr}
R_G^r
=
\sum_{\chi \in \mathrm{Irr}\, G} \FS_r (\chi)\, \chi,
\end{equation}
where $\mathrm{Irr}\, G$ denotes the set of (complex) irreducible characters of $G$ 
and $\FS_r (\chi)$ denotes the $r$-th Frobenius-Schur indicator of $\chi$.
Hence \eqref{InvSn} follows from  $\FS_2 (\chi) = 1$ 
$(\chi \in \mathrm{Irr}\, \Sn)$ and $\sum_{\chi \in \mathrm{Irr}\, \Sn} \chi (1) = \bll\STab (n)\brl$.

Let $H$ be a subgroup of $G$ and let $b$ be an element of $G$.
In this paper, we consider the following {\it coset-wise root number function}: 
$$ 
R_{G,bH}^r (a) 
:= \bll \bigl\{ c\in b H\,|\, c^r = a \bigr\} \bll.
$$
The support $K$ of the restricted function
$R_{G,bH}^r\bigl|_H$ becomes a subgroup of $H$ and the 
restriction of $R_{G,bH}^r$ on $K$ has the expansion
$$
R_{G,bH}^r \Bigl|_K
=
\sum_{\chi \in \mathrm{Irr} K} \FS_r (L_\chi) \chi,
$$
where $\FS_r (L_\chi)$ denotes the  
$r$-th Frobenius-Schur indicator of certain simple module $L_\chi$
of a weak Hopf algebra (WHA) $\FGX$ attached to $G$ and $X := G/H$. 
%When $r = 2$, we also show that $\nu_2 (L_\chi)$ agrees with a (generarized)
%twisted Frobenius-Schur indicator of Kawanaka-Matsuyama \cite{KM},
%unless $\nu_2 (L_\chi)\ne 0$.
When $G = \Sn$ and $H = \Sym_{n-1}$,
we give an explicit formulas of
$\FS_r (L_\chi)$ and $R_{G,bH}^r$ for every $r>0$.
When $G = \Sn$ and $H$ is a Young subgroup $\Sym_{\alpha}$,  we determine the value 
$\FS_2 (L_\chi)$ and give an explicit formula for the number of involutions in
$b \Sym_\alpha$.  
As a special case, we obtain
\begin{equation}
\label{bS}
\bll \bigl\{ a \in b \Sym_m \,|\, a^2 = 1 \bigr\} \bll
=
\begin{cases}
\,\bll\STab (m - k)\bll\,  \quad (\Sym_m b \Sym_m = \Sym_m b^{-1}\Sym_m)\\
\, 0  \qquad\qquad\qquad\quad (\Sym_m b \Sym_m \ne \Sym_m b^{-1}\Sym_m),
\end{cases}
\end{equation}
for each $0<m<n$ and $b \in \Sn$,
where 
$k := \bll \{1,2,\ldots, m\} \setminus b\{1,2,\ldots, m\}\bll$.
Also, we obtain
\begin{equation}
\label{bSS}
\bll \bigl\{ a \in b (\Sym_m \times \Sym_{m^\prime})\,|\, a^2 = 1 \bigr\} \bll
=
k !\, \bll\STab (m - k)\bll\, \bll\STab (m^\prime - k)\bll,
\end{equation}
where $m^\prime := n-m$. 

By counting $R^2_\Sn (1)$ using \eqref{bSS},  we obtain 
\begin{align}
\label{ST=sumST}
\bll\STab (n)\bll 
=
\sum_{0\leq k\leq \mathrm{min}\{m, m^{\prime} \}} 
\frac{m!\, m^{\prime} !}{k !\, (m-k)!\, (m^{\prime}-k)!}
\bll\STab (m - k)\bll\, \bll\STab (m^{\prime} - k)\bll. 
\end{align}
%

%Lusztig \cite{Lusztig} 

%\newpage

In \cite{NgSch}, Ng and Schauenburg defined Frobenius-Schur indicators as invariants of (objects of) pivotal fusion categories (See also \cite{Fuchs}). 
Also, Schauenburg \cite{Sch1,Sch2,Sch3}
gave several results for Frobenius-Schur indicators
of group-theoretical fusion categories.
Since the representation category of $\FGX$ is group-theoretical by 
Andruskiewitsch-Natale \cite{AN} and Mombelli-Natale \cite{MNatale}, 
our general results for Frobenius-Schur indicators of $\FGX$-modules overlap with
Schauenburg's results significantly. 
Nevertheless, we give WHA counterparts of his results, since WHA approach seems to be more elementary
than his category-theoretic approach.
In fact, our approach clarifies the importance of the G-sets $X$ and $X \times X$, 
which did not play important roles in his paper.

The outline of the paper is as follows. In Section 2 we give the definition of the algebra $\FGX$. 
In Section 3, we define and study Frobenius-Schur indicators of $\FGX$-modules.
In Section 4 we give a relation between second indicators of $\FGX$-modules and  
Kawanaka-Matsuyama indicators \cite{KM} of $\C K$-modules.
In Section 5 and Section 6, we compute second indicators of ${\cal F} (\Sn,\Sym_\alpha)$-modules
and indicators of ${\cal F} (\Sn,\Sym_{n-1})$-modules, respectively.
Also, we give the corresponding results for $R^r_{G,bH}$ in these sections.
In Section 7 we give a correspondence between invariant bilinear forms on $\FGX$-modules and
invariant bilinear pairings of some $\C K$-modules.  
In Section 8 we verify that our definition of Frobenius-Schur indicators of 
$\FGX$-modules coincides with that of Ng-Schauenburg \cite{NgSch}.

The author thank K. Shimizu for telling him about the Frobenius-Schur indicators.

\section{Preliminaries}
Throughout this paper, all modules are assumed to be finite dimensional over the
complex number field $\C$.
Let $G$ be a finite group and let $X$ be a finite left $G$-set.
For $x \in X$, we denote by $G_x$ the {\it stabilizer} of $G$ at $x$,
that is, $G_x = \{ a \in G\,|\, a x = x\}$. 
Let $\KXG$ be the $\C$-linear span of the symbols 
$e_x\, a\,\, (a \in G,\, x \in X)$.
Then $\KXG$ becomes an algebra via
$$
(e_x \, a)(e_y \, b) 
= 
\delta_{x,ay}\, e_x \, ab.
$$
By identifying $a \in G$ with $\sum_{x\in X} e_x\, a$, $\C  G$ becomes a subalgebra of 
$\KXG$.
The elements $e_x := e_x\, 1_G\,\, (x \in X)$ are mutually orthogonal idempotents and give a partition of unity of $\KXG$.

Let $M$ be a left $\KXG$-module and let $\Ome$ be an orbit of the $G$-set $X$.
We say that $M$ is of {\it type} $\Ome$ if $M = \bigoplus_{x \in \Ome} e_x M$.  
We note that each $\KXG$-module has a unique decomposition 
$M = \sum_{\Ome \in G\backslash X} M_\Ome$
such that $M_\Ome$ is of type $\Ome$.

It seems that the following is a folklore among some communities of Hopf algebraists.

\begin{prop}
\label{Ixy}
{\rm (cf.\,\cite{Lusztig} page 241, \cite{KMM} Lemma 3.2,\,Theorem 3.3)}
{\rm (1)}\, Let $\Ome = G x$ be an orbit of $X$ and let $V$ be a left $\C  G_{x}$-module.
Then $\II_{x} (V) := \C  G \otimes_{\C  G_{x}} V$ becomes a $\KXG$-modules via 
$$
a (b\otimes v) = ab \otimes v,\quad
e_y (b\otimes v) = \delta_{y,bx} b\otimes v\quad
(a,b \in G, y\in X, v \in V).
$$
{\rm (2)}\, The correspondence $\II_{x}$ gives an equivalence between 
the category of $\C  G_{x}$-modules
and the category of $\KXG$-modules of type $\Ome$.     
\end{prop}

Let ${\cal F} = \FGX$ be the $\C$-linear span of the symbols 
$e^x_y \, a\,\, (a \in G,\, x, y \in X)$.
Then ${\cal F}$ becomes an algebra via
$$
(e^x_y \, a)(e^z_w \, b)
=
\delta_{x,az} \delta_{y,aw} e^x_y \, ab. 
$$
Let 
$\Delta:\! {\cal F}\to {\cal F}\otimes {\cal F}$
and $\cu:\! {\cal F}\to \C $ be linear maps given by 
$$
\Delta (e^x_y \, a)
=
\sum_{z \in X} e^x_z \, a \otimes e^z_y \, a,
\qquad
\cu (e^x_y \, a) = \delta_{xy}.
$$
Then $\cal F$ becomes a $X$-face algebra with antipode
$S:\!{\cal F} \to {\cal F}; e^x_y a \mapsto a^{-1} e^y_x$ (cf. \cite{Hayashi}).
Hence $\cal F$ is a weak Hopf algebra (cf. \cite{WHAI}).
We call $\FGX$ the {\it group-like face algebra} of $(G,X)$. 

Let $\Ome = G(x,y)$ be an {\it orbital} of $X$, that is,
$\Ome \in G \backslash (X\times X)$.
Since $\FGX \cong \C ^{X\times X}\! \rtimes G$ as algebras,  
we have an equivalence $\II_{xy}$ between the category of $\C  G_{xy}$-modules 
and the category of $\FGX$-modules of type $\Ome$, where $G_{xy}$ stands for the 
{\it two-point stabilizer} $G_x \cap G_y$.
In particular, if
$\{ V(\lam) \}$ is a set of representatives for the isomorphism classes of simple
$\C  G_{xy}$-modules, then $\{ \II_{xy} (V(\lam))\}$ is a set of representatives for the isomorphism classes of simple %
$\FGX$-modules of type $\Omega$.

\section{Frobenius-Schur indicators}

We define elements $\inte^{[r]}\,\,(r\geq 1)$
of $\FGX$ by $\inte^{[1]} = \int := \frac{1}{|G|}\sum_{a\in G} \sum_{x \in X} e^x_x\, a$
and 
$\inte^{[r]} := (m^{(r)}\circ \Delta^{(r)}) (\inte)\,\, (r \geq 2) $
respectively,
where $m^{(r)}\!:\FGX^{\otimes r} \to \FGX$
and $\Delta^{(r)}\!:\FGX \to \FGX^{\otimes r}$ denote
the iterations of the product and the coproduct of $\FGX$ respectively, that is,
$m^{(3)} (\alp,\beta,\gam) =\alp \beta \gam$ and 
$\Delta^{(3)} (\alp) = (\Delta \otimes \id)(\Delta (\alp))$, for example.
Then, $\inte$ is an idempotent two-sided 
integral of 
$\FGX$ (cf.\,\cite{WHAI}), that is, $\inte^2 = \inte$ and
\begin{equation}
\label{defint}
\alpha \inte = \cu^L ( \alpha ) \inte,
\quad
\inte \alpha = \inte \cu^R (\alpha )
\quad
(\alpha \in \FGX),
\end{equation}
where, 
\begin{equation}
\label{eLdef}
\cu^L ( \alpha ) = \sum_{x,y,z\in X} \cu( e^x_z \alpha) e^z_y,
\quad
\cu^R ( \alpha ) = \sum_{x,y,z\in X} e^x_z \cu(\alpha e^z_y).
\end{equation}
Let $M$ be a finite-dimensional $\FGX$-module.
We define the $r$-{\it th Frobenius-Schur indicator} of $M$ by 
$\FS_r (M) := \tr_M (\inte^{[r]})$.
\begin{lem}
\label{mDS} 
{\rm (1)}\, Explicitly, the elements $\inte^{[r]}$ are given by
\begin{equation}
\label{sr=}
 \inte^{[r]} 
= 
\frac{1}{|G|}\sum_{a\in G}\sum_{x \in X} \delta_{x,a^r x}e^x_{a^{-1} x} a^r
=
\frac{1}{|G|}\sum_{a\in G}\sum_{x \in X} \delta_{x,a^r x} a^r e^x_{a^{-1} x}. 
\end{equation}
{\rm (2)}\, The element $\inte^{[r]}$
is central.  
\end{lem}

\par\noindent
{\it Proof.}\, Part (1) follows from the following computations:
\begin{align*}
& (m^{(r)}\circ \Delta^{(r)}) (e^x_x\, a)\\
= &
\sum_{y_1 \in X}\cdots \sum_{y_{r-1} \in X}
\left( e^x_{y_1} a \right) \left( e^{y_1}_{y_2} a \right) \cdots
\left( e^{y_{r-2}}_{y_{r-1}} a \right) \left( e^{y_{r-1}}_x a \right)\\
= &
\sum_{y_1 \in X}\cdots \sum_{y_{r-1} \in X}
\left( e^x_{y_1} a \right) \left( e^{y_1}_{y_2} a \right) \cdots
\left( e^{y_{r-3}}_{y_{r-2}} a \right) 
\delta_{y_{r-2},a y_{r-1}} \delta_{y_{r-1}, ax} 
\left( e^{y_{r-2}}_{y_{r-1}} a^2 \right)\\
= &
\sum_{y_1 \in X}\cdots \sum_{y_{r-3} \in X}
\left( e^x_{y_1} a \right) \left( e^{y_1}_{y_2} a \right) \cdots
\left( e^{y_{r-3}}_{a^2 x} a \right) 
\left( e^{a^2 x}_{ax} a^2 \right)\\
%= &
%\sum_{y_1 \in X}\cdots \sum_{y_{r-4} \in X}
%\left( e^x_{y_1} a \right) \left( e^{y_1}_{y_2} a \right) \cdots
%\left( e^{y_{r-4}}_{a^3 x} a \right) 
%\left( e^{a^3 x}_{a^2 x} a^3 \right)\\
= &\cdots\\
= &
\left( e^x_{a^{r-1} x} a \right) 
\left( e^{a^{r-1} x}_{a^{r-2} x} a^{r-1} \right) 
=
\delta_{x,a^r x}\,e^x_{a^{-1} x}\, a^r
=
\delta_{x,a^r x}\, a^r\, e^x_{a^{-1} x}.
\end{align*}
Let $c$ be an element of $G$. 
Replacing $a$ and $x$ in \eqref{sr=} by $cbc^{-1}$ and $cy$ respectively, 
we obtain
\begin{align*}
 \inte^{[r]} c
 = &
 \left(
 \frac{1}{|G|}\sum_{b\in G}\sum_{y \in X} \delta_{c y, c b^r y }\,
 e^{cy}_{c b^{-1} y}\, c b^r c^{-1}
 \right) c \\
 = &
 \frac{1}{|G|}\sum_{b\in G}\sum_{y \in X} \delta_{y,  b^r y }\,
 c\, e^{y}_{b^{-1} y}\, b^r 
 = 
 c \inte^{[r]}.
\end{align*}
For $y, z \in X$, we have
$$
 e^y_z \inte^{[r]} 
 =
 \frac{1}{|G|} \sum_{a\in G} \delta_{y,a^r y} \delta_{z,a^{-1}y}  e^y_{a^{-1} y} a^r
 =
 \frac{1}{|G|} \sum_{a\in G} \delta_{y,a^r y} \delta_{z,a^{-1}y}  a^r e^y_{a^{-1} y} 
 =
 \inte^{[r]} e^y_z.
$$ 
Since $g$'s and $e^y_z$'s generate the algebra $\FGX$, this proves Part (2). 
$\hfill \square$

\begin{theorem} {\rm (cf. Schauenburg \cite{Sch1}, Theorem 4.1)}
\label{FSIThm}
For each $x, y \in X$ and  $\C  G_{xy}$-module $V$,  we have 

\begin{equation}
\label{FSI}
\FS_r (\II_{xy} (V)) 
= 
\frac{1}{| G_{xy} |} \sum_{a \in G[x,y;\,r]} \tr_V (a^{-r}),
\end{equation}
where $G[x,y;r] := \{ a \in G\,|\, a x = y,\,  a^r x= x \}$. 
\end{theorem}

\par\noindent
{\it Proof.}
We first note that the right-hand side of \eqref{FSI} is well-defined, 
since $a^{-r} \in K:= G_{xy}$ for each $a \in G[x,y;\,r]$.
Also, we note that 
we may assume that $V$ is a simple $\C  K$-module,
since both the right-hand side and the left-hand side of \eqref{FSI}
are additive with respect to $V$.
Then, by Proposition \ref{Ixy} (2), 
$\II_{xy} (V)$ is a simple $\FGX$-module.
Hence by Schur's lemma, the action of the central element
$\inte^{[r]}$ 
on $\II_{xy} (V)$ 
is given by some scalar. Therefore, we have
\begin{gather*}
\tr_{\II_{xy} (V)} (\inte^{[r]}) 
= 
\frac{\dim \II_{xy} (V)}{\dim ( \C K \otimes_{\C K} V)}
\tr_{\C K \otimes_{\C K} V} (\inte^{[r]})\\
= 
\frac{|G|}{|K|}\tr_{\C K \otimes_{\C K} V} (\inte^{[r]}).
\end{gather*}
By Lemma \ref{mDS} (1), we have
\begin{gather*}
\inte^{[r]} (1_G \otimes v)
=
\frac{1}{|G|}
\sum_{a,z} \delta_{z,a^r z}\delta_{z, x}\delta_{a^{-1}z,y} a^r\otimes v \\
=
\frac{1}{|G|}
\sum_{c \in G[x,y;\,r]} 1_G \otimes c^{-r} v. 
\end{gather*}
Therefore, 
$$
\tr_{\C K \otimes_{\C K} V} (\inte^{[r]}) 
 = 
 \frac{1}{|G|}
\sum_{c \in G[x,y;\,r]} \tr_V (c^{-r}).
$$
This proves \eqref{FSI}.
$\hfill \square$

Let $H$ be a subgroup of $G$.
We define the {\it $r$-th root number function}
$R^r_{G}$ and the {\it $r$-th coset-wise root number function}
$R^r_{G,\,bH}$ by
\begin{gather*}
R^r_{G} (a) = \bll\{ c \in  G \,|\, c^r = a \} \bll,\\
R^r_{G,\,bH} (a) = \bll\{ c \in  b H \,|\, c^r = a \} \bll,
\end{gather*}
respectively.
We note that $R^r_G$ is a class function and that
\begin{equation}
\label{RrGhbH}
R^r_{G,\,hbH} (a) = R^r_{G,\,bH} (h^{-1} a h)
\end{equation}
for each $a,b \in G$ and $h \in H$.
In particular, we have
\begin{equation}
\label{R1=R1}
R^r_{G,\,hbH} (1) = R^r_{G,\,bH} (1).
\end{equation}
By \eqref{R1=R1}, the assignment 
$HbH \mapsto R^r_{G,HbH} (1):= R^r_{G,bH} (1)$
gives a well-defined function on $H\backslash G/H$.

\begin{prop}
The root number function satisfy
\begin{equation}
\label{R=sumR}
R^r_G (1) 
= 
\sum_{HbH \in H\backslash G/H} \frac{\bll H \bll}{\bll H \cap bHb^{-1} \bll}\,  
R^r_{G,\,HbH} (1).
\end{equation}

\end{prop}
\par\noindent
{\it Proof.}\, By definition, we have
\begin{align*}
R^r_G (1) 
& = 
\sum_{bH \in  G/H} R^r_{G,bH} (1)\\
& = 
\sum_{H b_1 H \in  H\backslash G/H} c_{H b_1 H} 
R^r_{G, H b_1 H} (1),
\end{align*}
where 
$c_{H b_1 H} := \bll \{bH \in  G/H\,|\, HbH = H b_1 H \}\bll$.
Since $c_{H b_1 H}$ is equals to the size of the $H$-orbit
through $y := b_1 H \in X:=G/H$,
it equals $|H|/|H_{y}|=|H|/|H \cap b_1 H b_1^{-1}|$.
This proves \eqref{R=sumR}.
$\hfill \square$

\begin{theorem}
{\rm (cf.\,Schauenburg \cite{Sch1}, Lemma 4.5)}
\label{RFmla}
For each $ a \in H$ and $y = b H  \in X := G/H$, we have 
\begin{gather}
\Bigl| \bigl\{ c \in  b H \,|\, c^r = a \bigr\} \Bigr|
\nonumber\\
\label{CRNum=}
=
\begin{cases}
 \sum_{\lam}\, \FS_r (\II_{xy} (V(\lam)))\, \chi_\lam (a) & (ay = y)\\
 0 & (ay \ne y), 
\end{cases}
\end{gather}
where ${x} = H \in X$,
$\{ { V(\lam)} \}$ is as in
Section 2 and $\chi_{\lam} = \tr_{V(\lam)}$ denotes the character of
$V(\lam)$.
\end{theorem}

\par\noindent
{\it Proof.}\,
To begin with, we show that the left-hand side of \eqref{CRNum=} is non-zero only if
$ay =y$.
Suppose that $c^r = a$ for some $c \in bH$.
Since $c \in bH$, we have $cx = b x =y$. 
Hence 
$$
ay = c^r y = c^{r+1} x = c a x = cx = y.
$$ 
Let $K$ be the two-point stabilizer $G_{xy}$.
By \eqref{RrGhbH}, $R^r_{G,bH}\bll_K$ is a class function on $K$.
Hence 
$$
R^r_{G,bH}\bll_K 
= 
\sum_\lam 
(R^r_{G,bH}\,|\,\chi_\lam )_K\, \chi_\lam, 
$$
where $(\,|\,)_K$ denotes the usual inner product of the space of class functions on $K$, that is, $(f|g)_K := |K|^{-1} \sum_{a\in K} f(a)\overline{g(a)}$. 
Therefore, it suffices to show that
$(R^r_{G,bH}\,|\,\chi_\lam )_K = \FS_r (\II_{xy} (V)(\lam)$.
By definition, we have
\begin{gather}
(R^r_{G,bH}\,|\,\chi_\lam )_K 
=
\frac{1}{|K|} \sum_{a \in K} 
\bll \{c \in bH\,|\, c^r = a\} \bll\,
\chi_\lam (a^{-1})
\nonumber\\
\label{(R|chi)}
=
\frac{1}{|K|} \sum_{c \in bH; c^r \in K} 
\chi_\lam (c^{-r}).
\end{gather}
Since
\begin{gather*}
\{ c \in bH\,|\, c^r \in K \}
=
\{ c \in G\,|\, cx =y, c^r y = y \}
=
G[x,y;r],
\end{gather*}
the right-hand side of \eqref{(R|chi)} coincides with that of 
\eqref{FSI} for $V = V(\lam)$ .
$\hfill \square$

\section{Twisted Frobenius-Schur indicators}

Let $K$ be a finite group. Let $\phi$ be an automorphism of $K$ and let $k_0$ be an element of $K$. We say that $(\phi, k_0)$ is an {\it outer involution} of $K$ 
if $\phi^2 (k) = k_0^{-1} k k_0\quad (k\in K)$ and $\phi (k_0) = k_0$. 
We note that if $K \leq G$ and $t \in G$ satisfies 
$t^{-1} K t = K$ and 
$t^2 \in K$, then 
$((-)^t, t^2)$ is an outer involution of $K$, where 
$(-)^t\!:K\to K; k \to t^{-1} k t$. 
Conversely, for an outer involution $(\phi,k_0)$,  
there exists a group $G \geq K$ and $t \in G \setminus K$ such that 
$(\phi, k_0) = ((-)^t, t^2)$.
Explicitly, $G$ is given by $G = K \coprod tK$, which is equipped with product 
$(tk)(t k^\prime) = k_0 (\phi (k)k^\prime)$, 
$(tk)k^\prime = t(kk^\prime)$, $k(tk^\prime) = t (\phi (k)k^\prime)$
$(k,k^\prime \in K)$, where $t K =\{ tk\,|\, k \in K\}$ is a copy of the set $K$. 

Let $V$ be a finite-dimensional $\C K$-module.
We define $(\phi,k_0)$-{\it twisted second Frobenius-Schur indicator} of $V$ by
$$
\FS_2 (V,\phi,k_0) = \frac{1}{|K|}\,\sum_{k \in K}\,\tr_V (k_0 \phi (k) k).
$$
When $(\phi, k_0) = ((-)^t, t^2)$, we write
$\FS_2 (V,t) = \FS_2 (V,\phi,k_0)$. 
It agrees with Kawanaka-Matsuyama's indicator (cf. \cite{KM}),
that is,
$$
\FS_2 (V,t) = \frac{1}{|K|}\,\sum_{k \in K}\,\tr_V ((tk)^2).
$$
We say that an orbital $\Ome$ is {\it symmetric} 
(or {\it self-paired}) if $\Ome^\tra = \Ome$, where
$\Ome^\tra := \{ (y,x)\,|\, (x,y) \in \Ome \}$.
\begin{prop}
{\rm (1)}
An orbital $\Ome =G (x,y)$ is symmetric if and only if there exists an element 
$t\in G$ such that $t x = y$, $ty = x$. 
In this case, $K:= G_{xy}$ satisfies 
$t^{-1} K t = K$ and $t^2 \in K$.\\
{\rm (2)} 
Let $H$ be a subgroup of $G$ and let $b$ be an element of $G$. 
Let $x_0$ be the element $H$ of $X = G/H$ and let $\Ome$ be $G (x_0,bx_0)$.
Then $\Ome$ is symmetric if and only if $HbH = Hb^{-1}H$.
\end{prop}
\par\noindent
{\it Proof.} Part (1) is obvious.
Since there exists a bijection
$H\backslash G /H \cong G\backslash (X\times X)$; $HbH \mapsto G(x_0,bx_0)$
(cf. \cite{CR} p240), $G(x_0,bx_0)$ is symmetric if and only if $HbH = Hb^{-1} H$. 
$\hfill \square$

\begin{theorem} {\rm (cf. Schauenburg \cite{Sch3}, Proposition 3.2)}
\label{FS2I}
Let $V$ be a $\C G_{xy}$-module.\\
{\rm (1)}\, If $\Ome = G(x, y)$ is not symmetric, then $\FS_2 (\II_{x,y} (V)) =0$.\\
{\rm (2)}\, Suppose that $\Ome$ is symmetric and that $t \in G$ satisfies
$t(x,y) = (y,x)$. Then, 
$$
\FS_2 (\II_{x,y} (V)) 
= 
\FS_2 (V,t).
$$ 
\end{theorem}

\par\noindent
{\it Proof.}
Since $G[x,y;2] = \{a \in G\,|\, ax = y, ay = x\}$ is empty if 
$\Ome^\tra \ne \Ome$, Part (1) follows from \eqref{FSI}. 
If $t \in G$ satisfies $t(x,y) = (y,x)$, then we have
$G[x,y;2] = G_{xy}\, t^{-1}$. Hence  Part (2) also follows from \eqref{FSI}. 
$\hfill \square$

\section{Symmetric groups I}

For each subset $S$ of $[n] := \{1,2,\ldots, n\}$,
we define a subgroup $\Sym (S)\cong \Sym_{|S|}$ of $G:=\Sn$ by
$\Sym (S):= \{ a\in \Sn \,|\, a i = i\,\, (i \in [n]\setminus S)\}$. 
For a set $S\subset \Z$ and an integer $\ep$, we set
$\ep + S = \{\ep + s\,|\, s\in S\}$.
Let $\alp = (\alp_1,\alp_2,\ldots,\alp_\ell)$ be a sequence of positive integers
such that $\alp_1 + \cdots + \alp_\ell = n$.
Let $\Sym_\alp = \Sym_{\alp_1}\times\cdots\times \Sym_{\alp_\ell}$ 
be the corresponding Young subgroup of $\Sn$.
Here, we identify $\Sym_\alp$ with  
$\Sym (A_1) \cdots \Sym (A_\ell)$ as usual,
where 
\begin{equation}
\label{defA}
A_1 = [\alp_1],\, A_2 = \alp_1+[\alp_2],\, \ldots,\,
A_\ell = \alp_1+\cdots \alp_{\ell-1} + [\alp_\ell].
\end{equation}
Next, define a set $X= \binom{[n]}{\alp}$ by 
$$
X
:=
\bigl\{ B = (B_1,\ldots,B_\ell) \in (2^{[n]})^\ell\,\bigm|\,
[n] = \coprod_i  B_i, \quad 
\bll B_i \bll = \alp_i\quad (1\leq i \leq \ell)\bigr\}.
$$
Then $X$ becomes a transitive $G$-set via
$a (B_1,\ldots,B_\ell) := (a(B_1),\ldots,a(B_\ell))$.
Since the stabilizer $G_A$ of $G$ at $A := (A_1,\ldots,A_\ell)$ is $\Sym_\alp$,
$\Sn/\Sym_\alp \cong X;$ $b\Sym_\alp \mapsto b A$ as $G$-sets. 
It is known that
$G (B,C) \mapsto [\bll B_i\cap C_j\bll]_{ij}$ gives a bijection 
from
$G \backslash(X \times X)$ onto 
$$
M_\alp
:=
\bigl\{ \Gam = [\gam_{ij}]_{ij}\in \mathrm{Mat}(\ell,\Z_{\geq 0})\, \bigm|\, 
\sum_i \gam_{ij} = \alp_j = \sum_i \gam_{ji}
\quad (1 \leq j \leq \ell)
\bigr\}.
$$
See e.g. \cite{Kerber}.
Note that $G (B, C)$ is a symmetric orbital if and only if 
$[\bll B_i\cap C_j\bll]_{ij}$ is a symmetric matrix.

Let $B = (B_1,\ldots,B_\ell) = b A$ be an element of $X$.
For $1 \leq i, j \leq \ell$,  we set 
$B_{ij} :=  A_i \cap B_j$ and
$\gam_{ij} := \bll B_{ij} \brl$. 
Also, we set $A_{ij} = \ep_{ij} + [\gam_{ij}]$, where 
$\ep_{11} =0$, $\ep_{12} =\gamma_{11}$,\,\ldots,
$\ep_{1\ell} =\gamma_{11}+\cdots + \gamma_{1,\ell-1}$,
$\ep_{21} = \gamma_{11}+\cdots + \gamma_{1\ell} = \alp_1$,
$\ep_{22} = \alp_1 + \gamma_{21}$,\,
\ldots,
$\ep_{ij} = \alp_{1}+\cdots + \alp_{i-1}+\gamma_{i1}+\cdots + \gamma_{i,j-1} $,\,\ldots\,.
By definition, we have $A_i = A_{i1}\coprod \cdots \coprod A_{i\ell}$ for each $1 \leq i\leq \ell$.
For each $i,j$, we fix a bijection $u_{ij}\!: A_{ij} \cong B_{ij}$
and define $u \in \Sn$ by $u\bll_{A_{ij}} = u_{ij}$.
Let
$\gam = (\gam_{1}, \gam_{2},\ldots, \gam_{\ell^2})$ be the sequence 
$(\gam_{11}, \gam_{12},\ldots, \gam_{1\ell}, \gam_{21},\ldots,\gam_{2\ell},\ldots,\gam_{\ell 1},\ldots,\gam_{\ell\ell})$
and let $K_0$ be the subgroup 
$\Sym_\gam = \Sym_{\gam_1}\times \cdots \times \Sym_{\gam_{\ell^2}}$ of 
$\Sn$, where $\Sym_0 = \{1\}$.
Note that we have
$a(\ep_{ij} + s) = \ep_{ij} + a_{ij}s$ for 
$a = (a_{11}, \ldots, a_{1\ell},\ldots, a_{\ell 1},\ldots,a_{\ell\ell})\in \Sym_\gam$ and
$s \in [\gam_{ij}]$.
Let 
$K$ the two-point stabilizer $G_{AB}$.

\begin{lem}
{\rm (1)}\, $K = \prod_{ij} \Sym (B_{ij})$.\\
{\rm (2)}\, The correspondence $k \mapsto u^{-1} k u$ gives a group isomorphism
$\psi\!:K \cong K_0$.
\end{lem}
{\it Proof.} (1)
For each $a \in \Sn$, $a \in K$ if and only if $a A_i = A_i$ and $a B_j = B_j$ for all $i,j$, if and only if $a B_{ij} = B_{ij}$ for all $i,j$.
Hence, we have Part (1).\\
(2)  Since $K_0 = \prod_{ij}\Sym (A_{ij})$ and $\bll \Sym (A_{ij})\bll = \gam_{ij}! = \bll \Sym (B_{ij})\bll$,
it suffices to show that $u k u^{-1} \in \Sym (B_{ij})$ for each $k \in \Sym (A_{ij})$. 
Let $s$ be an element of $[n] \setminus B_{ij}$. 
Since $u^{-1} s \in [n]\setminus A_{ij}$, we have $k u^{-1} s = u^{-1} s$.
This proves the assertion.
$\hfill \square$

\bigskip
\par\noindent
Now suppose that $G(A,B)$ is a symmetric orbital.
We define  $t_0, t \in \Sn$ by 
$t_0 (\ep_{ij} + s) = \ep_{ji} + s$\quad $(s \in [\gam_{ij}])$ and
$t = u t_0 u^{-1}$.
Then we have $t_0^2 = 1_G$, $t_0 \bll_{A_{ii}} = \id$ and
$t_0 (A_{ij}) = A_{ji}$.
Moreover, 
$t_0 a t_0 = a^\tra$ for 
$a =  (a_{11}, \ldots, a_{1\ell},\ldots, a_{\ell 1},\ldots,a_{\ell\ell})\in K_0$,
where 
$$a^\tra := (a_{11}, a_{21},\ldots, a_{\ell1}, ,\ldots,a_{1\ell},\ldots,a_{\ell\ell}).
$$
Since $t(B_{ij}) = B_{ji}$, we have $t A = B$, $t B = A$.
Moreover, we have $t \psi^{-1} (a) t = \psi^{-1} (a^\tra)$ for each $a \in K_0$.

For $m \geq 0$, let $\PP (m)$ be the set of partitions of $m$ and let 
$\{ V(\lam)\,|\, \lam \in \PP (m)\}$ be a complete representatives of simple
$\C \Sym_m$-modules such that $\dim V(\lam) = \bll \STab (\lam)\bll$,
where $\STab (\lam)$ denotes the set of standard tableaux of shape $\lam$.  
We denote the character of $V(\lam)$ by $\chi_\lam$. 
Note that $\PP (0)$ is a single element set $\{ ()\}$
and that $V(())$ is a one-dimensional module of $\Sym_0 = \{ 1 \}$.
 
Let 
$\PP (\Gam)$ be the set of matrices
$\Lam = [\lam_{ij}]_{1\leq i,j \leq \ell}$ of partitions such that 
$\lam_{ij} \in \PP (\gam_{ij})$. 
For each $\Lam =[\lam_{ij}] \in \PP (\Gam)$, 
define a simple $\C K_0$-module $V(\Lam)$ by the following outer tensor product: 
$$
V(\Lam) = V(\lam_{11})\boxtimes V(\lam_{12})\boxtimes\cdots\boxtimes V(\lam_{1\ell})\boxtimes \cdots \boxtimes V(\lam_{\ell 1})\boxtimes\cdots\boxtimes V(\lam_{\ell\ell}).
$$
Then $\{ V(\Lam)\,|\, \Lam \in \PP (\Gam)\}$ gives a complete representative of simple
$\C K_0$-modules. Hence 
$\{ V(\Lam)^{\psi}\,|\, \Lam \in \PP (\Gam)\}$ 
and
$\{ {\cal I}_{A,B}(V(\Lam)^{\psi})\,|\, \Lam \in \PP (\Gam)\}$
give complete representatives of simple
$\C K$-modules and simple $\FGX$-modules of type $G(A,B)$, respectively.
Here, the action of $K$ on $V(\Lam)^{\psi}:=V(\Lam)$ is given by
$(k,v) \mapsto \psi (k) v$\,\, ($k \in K$, $v \in V(\Lam)^{\psi}$).

\begin{theorem}
\label{FS2YSubG}
For each $\Lam \in \PP (\Gam)$, the second Frobenius-Schur indicator of
the $\FGX$-module 
${\cal I}_{A,B}(V(\Lam)^{\psi})$
is given by
\begin{equation}
\mathrm{FS}_2 ({\cal I}_{A,B} (V(\Lam)^\psi))
=
\begin{cases}
 1\, & (\Gam^\tra = \Gam\,\, \mathrm{and}\,\, \Lam^\tra =\Lam)\\
 0 & (\Gam^\tra \ne \Gam\,\,\,\, \mathrm{or}\,\,\,\, \Lam^\tra \ne \Lam).
\end{cases}
\end{equation}
\end{theorem}
\bigskip
{\it Proof.}
By Theorem \ref{FS2I} (1), we may assume $\Gam^\tra = \Gam$.
Hence 
$$
\mathrm{FS}_2 ({\cal I}_{A,B} (V(\Lambda)^\psi))
= \mathrm{FS}_2 (V(\Lambda)^\psi,t)
= \frac{1}{\bll K \bll} \sum_{k\in K}
\tr_{V(\Lam)} (\psi ((tk)^2))\\
$$
by Theorem \ref{FS2I} (2).
Since $\psi ((tk)^2) = t_0 a t_0 a = a^{\!\tra} a$ for $k = \psi^{-1} (a)$
and $a \in K_0$,
the right-hand side equals
\begin{align*}
& 
\frac{1}{\bll K_0 \bll} \sum_{(a_{11},\ldots,a_{\ell\ell})\in K_0}
\tr_{V(\Lam)} (a_{11} a_{11},a_{21}a_{12},\ldots,a_{\ell 1}a_{1 \ell},
\ldots a_{1 \ell}a_{\ell 1},\ldots, a_{\ell\ell}a_{\ell\ell})\\
= & 
\frac{1}{\prod_{ij}\gam_{ij}!} 
\sum_{a_{11}\in \Sym_{\gam_{11}}}
\sum_{a_{12}\in \Sym_{\gam_{12}}}\cdots
\sum_{a_{\ell\ell}\in \Sym_{\gam_{\ell\ell}}}
\prod_{ij}\chi_{\lam_{ij}} (a_{ji} a_{ij}).
\end{align*}
Since 
$\frac{1}{\gam_{ii}!} \sum_{a_{ii}} \chi_{\lam_{ii}} (a_{ii}^2) = \FS_2 (V(\lam_{ii})) = 1$,
the right-hand side equals
\begin{align*}
& \prod_{i < j}
\frac{1}{(\gam_{ij}!)^2} 
\sum_{a,a^{\prime}\in \Sym_{\gam_{ij}}}
\chi_{\lam_{ij}} (a^{\prime} a)
\chi_{\lam_{ji}} (a a^{\prime})\\
= &
\prod_{i < j}
\frac{1}{\gam_{ij}!} 
\sum_{a\in \Sym_{\gam_{ij}}}
\chi_{\lam_{ij}} (a )
\overline{\chi_{\lam_{ji}} (a)}\\
= & 
\prod_{i < j}
\delta_{\lam_{ij},\lam_{ji}}
=
\delta_{\Lam,\Lam^\tra},
\end{align*}
where the first equality follows from the fact that
$\chi_{\lam_{ji}}$ is a real-valued class function.
The second equality follows from 
the orthogonality relation of irreducible characters.
$\hfill \square$

\begin{theorem}
Let $\alp =(\alp_1,\ldots,\alp_\ell)$ be a sequence of positive integers such that 
$\sum_i \alp_i = n$ and let $\Sym_\alp$ be the corresponding Young subgroup of $\Sn$.
Then, for each $b \in \Sn$, we have
\begin{equation}
\label{FS2Y}
\bigl| \bigl\{ a \in  b\, \Sym_{\alpha}\,|\, a^2 = 1 \bigr\} \bigr|
=
\begin{cases}
 \left(\prod_{i < j} \gamma_{ij} !\right)\,
 \prod_{i} \bigl| \mathrm{STab} (\gamma_{ii}) \bigr|
  & (\Gam^\tra = \Gam)\\
 0 & (\Gam^\tra \ne \Gam), 
\end{cases}
\end{equation}
where $\Gam = [\gam_{ij}]$, $\gam_{ij} = \bll A_i \cap b A_j \bll$ and $A_i$ is as in
\eqref{defA}.
\end{theorem}

\par\noindent
{\it Proof.} By Theorem \ref{RFmla} and Theorem \ref{FS2YSubG}, the left-hand side
of \eqref{FS2Y} equals
\begin{align*}
& 
\sum_{\Lam = [\lam_{ij}]\in \PP (\Gam)%;\,\, (1\leq i,j\leq \ell);\Lam^\tra=\Lam
}
\delta_{\Gam, \Gam^\tra}\, \delta_{\Lam, \Lam^\tra} \dim V(\lam_{11}) \dim V(\lam_{12})
\cdots \dim V(\lam_{\ell\ell})\\
= & \,\,\delta_{\Gam, \Gam^\tra} \Bigl( \prod_{i} 
\sum_{\lam_{ii} \in \PP (\gam_{ii}) } \bll\STab(\lam_{ii}) \bll \Bigr)
\,\,
\Bigl( \prod_{i < j} 
\sum_{\lam_{ij} \in \PP (\gam_{ij}) } 
\bll\STab(\lam_{ij}) \bll^2 \Bigr)\\
= & \,\,\delta_{\Gam, \Gam^\tra} \Bigl(\prod_{i}\, \bll\STab(\gam_{ii}) \bll \Bigr)
\,\,
\Bigl( \prod_{i < j} 
\gam_{ij} ! \Bigr)
\end{align*}
as desired.   
$\hfill \square$

\par\noindent
{\it Example.}
(1) Suppose $\alp = (m,m^\prime)$, where $m^\prime = n-m$.
Since $\Gam = [\gam_{ij}] \in M_\alp$ satisfies 
$\gam_{11} + \gam_{12} = m = \gam_{11} + \gam_{21}$, 
it is a symmetric matrix of the form
$$
\begin{bmatrix}
m-k &k\\
k & m^\prime - k
\end{bmatrix}.
$$
Hence each orbital of $X \cong \Sn/ (\Sym_m\times \Sym_{m^\prime})$ is symmetric.
When $\Gam$ corresponds to $G(A,bA)$, $k$ equals 
$\bll A_1 \cap bA_2\bll =\bll [m]\setminus b[m] \bll$.
Hence the number of involutions in $b(\Sym_m\times \Sym_{m^\prime})$ is given by
\eqref{bSS}. 
Since $\bll H \cap b H b^{-1}\bll = \bll G_{A,bA} \bll = (m-k)! (k!)^2 (m^\prime-k)!$, 
\eqref{ST=sumST} follows from 
\eqref{R=sumR}.
\\
(2) Suppose $\alp = (m,1^{n-m})$, so that $\Sym_\alp \cong \Sym_m$
and $A = ([m],\{ m+1 \},\ldots,\{ n \})$. 
Then, $\gam_{ij} = 0,1$ unless $(i,j) = (1,1)$.
Assume $\Gam$ is symmetric. Since $\bll \PP (0) \bll = \bll \PP (1) \bll = 1$,
every matrix $\Lam \in \PP (\Gam)$ is necessarily symmetric.
Hence $\FS_2 ({\cal I}_{A,B} (V)) = 1$ for each simple $G_{AB}$-module $V$.
We note that this result gives a characterization of ``null indicator double coset" of 
Schauenburg \cite{Sch3} Theorem 4.2.
Also, we note that this gives a generalization of computations of Frobenius-Schur indicators of
Hopf algebra representations given by 
Kashina, G. Mason, S. Montgomery \cite{KMM},
Jedwab, S. Montgomery \cite{JM} and Timmer \cite{Timmer} (cf. \cite{Sch3}).
Since $\gam_{11} = m - \bll [m]\setminus b[m]\bll$, we obtain \eqref{bS}.

\section{Symmetric groups II}

Let $G = \Sn$ be the symmetric group of $n$-letters. 
As usual, we identify the two-point stabilizer $G_{n,n-1}$ with
$\Sym_{n-2}$.

\begin{theorem}
For each $\C \Sym_{n-2}$-module $V$, 

\begin{equation}
\label{FSrI}
\FS_r (\II_{n,n-1}(V))
=
\sum_{2 \leq s \leq n;\,\,s|r} 
\FS_r ( V|_{\mathfrak{S}_{n-s}}),
\end{equation}
where $V|_{\mathfrak{S}_{n-s}}$ denotes the restriction of $V$
to $\Sym_{n-s}$.
Here, for convenience, we set $\FS_r ( V|_{\mathfrak{S}_{0}}) = \dim V$.

\end{theorem}

\par\noindent
{\it Proof.}\,
Let $a$ be an element of $G[x,y;r]$,
where $G[x,y;r]$ is as in Theorem \ref{FSIThm},
$i_1 = x: = n$ and $i_2 = y := n-1$.
Let $s \geq 2$
be the smallest integer such that $a^s x = x$.
It is easy to see that $s$ divides $r$ and that
$s$ agrees with the size of the orbit $\langle a \rangle  x$.
Hence we have the following decomposition:
\begin{gather}
G[x,y;r] 
=
\coprod_{2 \leq s \leq n;\,\,s|r} 
G_s [x,y;r], \\
G_s [x,y;r] := \bigl\{ a \in G[x,y;r]\,\,\bll\,\, 
\bll \langle a \rangle  x\bll = s \bigr\}. 
\end{gather} 
Suppose that $a$ belongs to $G_s [x,y;r]$.
We define integers $i_3,i_4,\ldots,i_{s} \in [n-2]$ by
$i_{3} = a^{2} x,\ldots,i_s = a^{s-1} x$.
Then, we have 
$h:=a(i_1,i_2,\ldots,i_{s})^{-1} \in \Sym([n]\setminus \{i_1,\ldots,i_s\})$,
that is, $h$ fixes each element of $\{i_1,\ldots,i_s \}$. 
Conversely, if $i_3,\ldots,i_{s}$ are distinct elements of $[n-2]$ and 
$h$ is an element of $\Sym([n]\setminus \{i_1,\ldots,i_s\})$,
$a = h (i_1,\ldots,i_{s})$ gives an element of $G_s [x,y;r]$.
Therefore, we have
\begin{gather}
G_s [x,y;r] 
=
\coprod_{i_3,\ldots,i_{s}} 
\bigl\{ h (i_1,\ldots,i_{s}) \,\,\bll\,\, 
h \in \Sym([n]\setminus \{i_1,\ldots,i_s\}) \bigr\}, 
\end{gather} 
where $i_3,\ldots,i_{s}$ run over distinct elements of $[n-2]$.
By \eqref{FSI}, this implies
\begin{align*}
& \FS_r (\II_{n,n-1}(V))\\
=\,\, &\frac{1}{(n-2)!}\sum_{2 \leq s \leq n;\,\,s|r}\, 
\sum_{i_3,\ldots,i_{s}} \sum_{h}
\tr_V ((h (i_1,\ldots,i_{s}))^{-r}).
\end{align*}
where $h$ runs over $\Sym([n]\setminus \{i_1,\ldots,i_s\})$.
Since $h$ commutes with the permutation
$(i_1,\ldots,i_{s}) \in \Sym(\{i_1,\ldots,i_s \})$,
the right-hand side of the above equality is
\begin{gather*}
\frac{1}{(n-2)!}\sum_{2 \leq s \leq n;\,\,s|r}\, 
\sum_{i_3,\ldots,i_{s}} \sum_{h}
\tr_V (h^{-r})\\
= \frac{1}{(n-2)!}\sum_{2 \leq s \leq n;\,\,s|r} 
\frac{(n-2)!}{(n-s)!} \sum_{h\in \Sym_{n-s}}
\tr_V (h^{-r})\\
= \sum_{2 \leq s \leq n;\,\,s|r} 
\FS_r ( V|_{\mathfrak{S}_{n-s}}). 
\end{gather*}
$\hfill \square$
\begin{theorem}
{\rm (1)}\, 
We have
\begin{equation}
\label{R=sumInd}
R^r_{\Sn,(n-1,n)\Sym_{n-1}}|_{\mathfrak{S}_{n-2}}
=
\sum_{2 \leq s \leq n;\,\,s|r} 
  \Ind_{\mathfrak{S}_{n-s}}^{\mathfrak{S}_{n-2}}(R^r_{\mathfrak{S}_{n-s}}),
\end{equation}
where $\Ind_{\mathfrak{S}_{n-s}}^{\mathfrak{S}_{n-2}}(R^r_{\mathfrak{S}_{n-s}})$
denotes the induced class function of $R^r_{\mathfrak{S}_{n-s}}$
on $\Sym_{n-2}$.\\
{\rm (2)}\, 
The class function $R^r_{\Sn,(n-1,n)\Sym_{n-1}}|_{\mathfrak{S}_{n-2}}$
is a character of a certain representation of $\mathfrak{S}_{n-2}$.
\end{theorem}

\par\noindent
{\it Proof.}
Let $H$ be a finite group.
By \eqref{RGr} and the orthogonal relation of the irreducible characters, 
we have 
$\FS_r (W) = (R^r_H | \chi_W)_H$
for each simple $\C H$-module $W$.
Since $\FS_r$ is additive, this relation also holds for every finite-dimensional
$\C H$-module $W$. 
Suppose $H$ is a subgroup of a finite group $G$. By Frobenius reciprocity, 
we have
$\FS_r(V \bll_H) = (\Ind^G_H (R^r_H)\,|\,\chi_V)_G$
for each finite-dimensional $\C G$-module $V$. 
Applying this equality to $G = \Sym_{n-2}$ and $H = \Sym_{n-s}$ together with \eqref{FSrI},
we find that
$$
\FS_r (\II_{n,n-1}(V))
=
\sum_{2 \leq s \leq n;\,\,s|r} 
(\Ind_{\mathfrak{S}_{n-s}}^{\mathfrak{S}_{n-2}}(R^r_{\mathfrak{S}_{n-s}})\,|\,\chi_V)_{\Sym_{n-2}}.
$$
Hence by Theorem \ref{RFmla}, we get 
\begin{align*}
R^r_{\Sn,(n-1,n)\Sym_{n-1}}|_{\mathfrak{S}_{n-2}}
& =\sum_{\mu \in \PP_{n-2}} 
\sum_{2 \leq s \leq n;\,\,s|r} 
(\Ind_{\mathfrak{S}_{n-s}}^{\mathfrak{S}_{n-2}}(R^r_{\mathfrak{S}_{n-s}})\,|\,\chi_\mu)_{\Sym_{n-2}}\, \chi_\mu\\
& =
\sum_{2 \leq s \leq n;\,\,s|r} 
  \Ind_{\mathfrak{S}_{n-s}}^{\mathfrak{S}_{n-2}}(R^r_{\mathfrak{S}_{n-s}}),
\end{align*}
where the last equality follows from the fact that 
$\{ \chi_\mu\,|\,\mu \in \PP_{n-2}\}$ is an orthonormal basis of the space of class functions of $\Sym_{n-2}$.
Thus we get Part (1). 
Part (2) follows immediately from Part (1) and  \cite{SnInd}.
$\hfill \square$

\begin{cor}
{\rm (1)}\, 
For each $b \in \Sn \setminus \Sym_{n-1}$, 
%
%$$
%R^r_{\Sn, b\Sym_{n-1}} (1) 
%= 
%\sum_{2\leq s \leq n;\, s | r}
%\frac{(n-2)!}{(n-s)!} R^r_{\Sym_{n-s}} (1).
%$$
%
\begin{equation}
\label{RNFmla}
\bll \bigl\{ a \in b \Sym_{n-1}\,|\, a^r = 1 \bigr\} \bll
=
\sum_{2\leq s \leq n;\, s | r}
\frac{(n-2)!}{(n-s)!}\,
\bll \bigl\{ a \in  \Sym_{n-s}\,|\, a^r = 1 \bigr\} \bll.
\end{equation}
{\rm (2)}\, {\rm (cf. \cite{Chowla})}
The root number $R^r_{\Sn} (1)$ satisfies the recurrence relation
\begin{equation}
\label{Recurs}
R^r_{\Sn} (1) =
\sum_{1\leq s \leq n;\, s | r}
\frac{(n-1)!}{(n-s)!}\,
R^r_{\Sym_{n-s}} (1).
\end{equation}
\end{cor}

\par\noindent
{\it Proof.}
Since the induced class function $\Ind^G_H (f)$ satisfies
$\Ind^G_H (f)(1) = \frac{|G|}{|H|} f(1)$, 
\eqref{RNFmla} immediately follows from  \eqref{R=sumInd} when
$b = (n-1,n)$.
On the other hand, since  
$\Sym_{n-1} \backslash \Sn / \Sym_{n-1} 
= \{ \Sym_{n-1}, \Sym_{n-1}(n-1,n)\Sym_{n-1}\}$,  
we have $R^r_{G, b\Sym_{n-1}} (1) = R^r_{G,(n-1,n) \Sym_{n-1}} (1)$
by \eqref{R1=R1}. This proves Part (1).
Part (2) follows from Part (1) and \eqref{R=sumR}.
$\hfill \square$

%\newpage

\section{A correspondence between bilinear pairings}

As well as the group case, the second Frobenius-Schur indicator of an $\FGX$-module 
$M$ has a close connection to invariant bilinear forms on $M$.  
In this section, we show it by giving a correspondence between 
invariant bilinear forms 
on $\FGX$-modules and certain bilinear pairings on $\C G_{xy}$-modules. 

\bigskip
Let $G$ be a finite group and let $K$ be its subgroup.
Let $t$ be an element of $G$ such that 
$t^{-1} K t = K$ and $t^2 \in K$.
For a $\C K$-module $V$, we denote by 
${}^t V =\{ {}^t v\,|\,v\in V\}$ a copy of $V$ with $\C  K$-action given by
$k\, {}^t v :=  {}^t (t^{-1} k t v)$.
Let $B : V \times {}^t V \to \C $
be a bilinear pairing. We say that $B$ is $K$-{\it invariant} if 
$B(k v,{}^t w) = B(v, k^{-1}\, {}^t w)$ for each $k \in K$ and $v,w \in V$.
We denote by ${\cal B}(V,t)$ the set of $K$-invariant bilinear pairings  
$B : V \times {}^t V \to \C $.
For $B \in {\cal B}(V,t)$, we set $B^\tra (v, {}^t w) := B(t^2 w, {}^t v)$. 
Since
\begin{gather*}
B^\tra (kv, {}^t w) = B(t^2 w, t k t^{-1}\, {}^t v)
= B(t k^{-1} t w,  {}^t v)
= B^\tra (v, k^{-1}\, {}^t w),\\
(B^\tra)^\tra (v, {}^t w) = B(t^2v, t^2\, {}^t w) = B (v,{}^t w),
\end{gather*}
we have $B^\tra \in {\cal B}(V,t)$ and 
$(B^\tra)^\tra = B$.
Similarly to \cite{KM}, we have the following result.

\begin{prop}
\label{FS2V}
Let $V$ be a simple $\C K$-module.
Then $\FS_2 (V,t) \in \{ 0, \pm 1 \}$ and $\dim {\cal B} (V,t) \leq 1$.
Moreover,  we have
\begin{equation}
 \FS_2 (V,t) =
 \begin{cases}
 1  & \dim {\cal B} (V,t)_+ = 1\\
 -1 & \dim {\cal B} (V,t)_- = 1\\
 0  & \dim {\cal B} (V,t) = 0,\\
 \end{cases}
\end{equation}
where ${\cal B} (V,t)_\pm := \{ B \in {\cal B} (V,t)\,|\, B^\tra = \pm B \}$. 
\end{prop}

Let $M$ be a $\FGX$-module and let $C\!: M \times M \to \C $ be a 
bilinear form on $M$. We say that $C$ is $\FGX$-{\it invariant} if it is $G$-invariant and satisfies
$C( e^x_y \xi, \eta ) = C(\xi, e^y_x \eta)$ for each $x, y \in X$ and $\xi, \eta \in M$.
We denote by  ${\cal B}(M)$  the set of $\FGX$-invariant bilinear forms.
For $C \in  {\cal B}(M)$, we define $C^\tra \in {\cal B}(M)$ by
$C^\tra (\xi, \eta ) = C(\eta,\xi)$.

\begin{thm}
\label{BCorresp}
Let $\Omega = G(x,y)$ be a symmetric orbital and let 
$t$ be an element of $G$ such that $t(x,y) = (y,x)$.\\
{\rm (1)}\, For each $\C G_{xy}$-module
$V$, we have a bijective correspondence
$$
\Res\!: {\cal B} (\II_{x,y} (V)) \cong {\cal B}(V,t)
$$
given by
$\Res (C)(v,{}^t w) = C(1\otimes v, t\otimes w)$.
The inverse $\Ind$ of $\Res$ is given by
$\Ind (B) (a\otimes v, b \otimes w) 
= 
\sum_{k \in K} \delta_{akt, b}\, B(v, k\,  {}^t w)$,
where $K = G_{xy}$.\\
{\rm (2)}\, For each $C \in {\cal B} (\II_{x,y} (V))$, we have
$\Res(C^\tra) = \Res(C)^\tra$.\\
{\rm (3)}\, A pairing $B \in {\cal B} (V,t)$ is 
non-degenerate if and only if $\Ind (B)$ is non-degenerate. 
\end{thm}
\par\noindent
{\it Proof.}\, It is straightforward to verify that   
$\Res$ and $\Ind$ give well-defined maps between 
${\cal B} (\II_{x,y} (V))$ and ${\cal B}(V,t)$.
Also, it is easy to verify that $\Res \circ \Ind = \id$.
Hence, to show Part (1), it suffices to prove that 
\begin{equation}
\label{IndRes=Id}
\Ind(\Res (C)) (a\otimes v, b \otimes w) = C(a\otimes v, b \otimes w)
\end{equation}   
for each $a, b \in G$ and $v, w \in V$.
By $\FGX$-invariance, the left and right-hand sides of \eqref{IndRes=Id}
is zero unless $a(x,y) = b(y,x)$.
Suppose $a(x,y) = b(y,x)$.
Since $k:=a^{-1} b t^{-1} \in K$, the left-hand side of \eqref{IndRes=Id}
is
$$
\Res (C)(v,k\,{}^t w) 
=
C (1 \otimes v, t \otimes t^{-1} k t w)
= 
C (1 \otimes v, k t \otimes w)
=
C(a\otimes v, b \otimes w).
$$ 
This proves Part (1).
Part (2) follows from $t^2 \in K$.
Suppose that $C = \Ind (B) \in {\cal B} (\II_{x,y} (V))$ is non-degenerate and that
$v \in V$ satisfies $B(v,{}^t w) = 0$ for every $w \in V$.
To prove the non-degeneracy of $B$,
it suffices to show that 
$C(1 \otimes v, b \otimes w) = 0$ for each $b \in G$ and
$w \in V$. 
By $\FGX$-invariance, we may assume $b (x, y) = (y,x)$,
or $b = t k$ for some $k \in K$.
Then, we obtain   
$C(1 \otimes v, b \otimes w) = B (v,{}^t (kw)) = 0$
and prove the non-degeneracy of $B$.
Conversely, suppose that $B$ is non-degenerate and that 
$m \in e^z_w \II_{x,y} (V)$ satisfies
$C(m, n) = 0$ for every $n \in e^w_z \II_{x,y} (V)$. 
Let $a$ be an arbitrary element of $G$ such that $a(x,y) = (z,w)$. 
Then, we have that $m = a \otimes v$ and $n = at \otimes w$ for some $v, w \in V$,
and that
$B(v,{}^t w) = C(m, n) = 0$.
Hence, the non-degeneracy of $B$ implies $m = 0$.
Since $C$ is $\FGX$-invariant, this proves the non-degeneracy of $C$.
$\hfill \square$

\begin{cor}
\label{FS2M}
Let $M$ be a simple $\FGX$-module.
Then $\FS_2 (M) \in \{ 0, \pm 1 \}$ and $\dim {\cal B} (M) \leq 1$.
Moreover,  we have
\begin{equation}
 \FS_2 (M) =
 \begin{cases}
 1  & \dim {\cal B} (M)_+ = 1\\
 -1 & \dim {\cal B} (M)_- = 1\\
 0  & \dim {\cal B} (M) = 0,\\
 \end{cases}
\end{equation}
where ${\cal B} (M)_\pm := \{ C \in {\cal B} (M)\,|\, C^\tra = \pm C \}$. 
\end{cor}

\par\noindent
{\it Proof.}\, Suppose $M$ is of type $\Omega$. When $\Omega$ is symmetric,
the assertion follows immediately from 
Theorem \ref{FS2I} (2), Proposition \ref{FS2V} and Theorem \ref{BCorresp}.
When $\Omega$ is not symmetric,
${\cal B} (M) = 0$ by the definition of $\FGX$-invariance.
Hence the assertion follows from Theorem \ref{FS2I} (1).  
$\hfill \square$

\begin{prop}
Let $x, y $ and $b$ be as in Theorem {\rm \ref{RFmla}} and let
$\Ome$ be $G(x,y)$. Then the following conditions are equivalent: \\
{\rm (1)}\, $\Ome^\tra = \Ome$.\\ 
{\rm (2)}\, $R^2_{G,bH}\ne 0$.\\
{\rm (3)}\, $\FS_2 (M) \ne 0$ for some $\FGX$-module $M$ of type $\Omega$.  
\end{prop}

\par\noindent
{\it Proof.}\, The equivalence of (2) and (3) follows from 
Theorem \ref{RFmla} and the linear independence of the characters.
Since the unit $\C G_{xy}$-module $\C$ satisfies
$\dim {\cal B} (\C,t)_+ = 1$,
the equivalence of (1) and (3) follows from Theorem \ref{BCorresp} (1),
Corollary \ref{FS2M} and Theorem \ref{FS2I} (1).
$\hfill \square$

\section{Frobenius-Schur indicators of Ng and Schauenburg}

In \cite{NgSch}, Ng and Schauenburg have defined higher Frobenius-Schur indicators
$\nu_r (M)$
for each pivotal tensor category ${\cal C}$ and its object $M$.
In this section, we verify that $\FS_r$ coincides with $\nu_r$ when
${\cal C}$ is the category ${}_{\FFF} {\bf mod}$
of finite-dimensional left $\FFF$-modules, where $\FFF =\FGX$ for some $(G,X)$.
We refer to \cite{EGNO} for terminology for tensor categories.
To begin with, we give an explicit description of operations on ${}_{\FFF} {\bf mod}$.
For each $M, N \in \mathrm{ob} {}_{\FFF} {\bf mod}$,  
let $M \tb N$ be a subspace  $M \otimes N$ defined by
$M \tb N := 
\Delta (1) (M\otimes N) = \bigoplus_{z\in X} e_z M \otimes e^z N$,
where $e_y = \sum_x e^x_y$ and $e^x = \sum_y e^x_y$. 
Then $M \tb N$ becomes an $\FFF$-module via
$$
e^x_y\,a \sum_z e_z m \otimes e^z n 
=  \sum_z e_z e^x a m \otimes e^z e_y a n
\,\,\, (a \in G,\, x, y,z \in X,\,
m \in M, n \in N). 
$$
The linear span $\uob := \C X$ of $X$ becomes an $\FFF$-module via
$e^x_y\, a \otimes z \mapsto \delta_{x,az}\delta_{y,az}az$.
Moreover it becomes a unit object with respect to $\tb$ via
\begin{align}
M \cong \uob \tb M 
& = 
\bigoplus_x \C x \otimes e^x M;\,\,m \mapsto \sum_x x \otimes e^x m,
\nonumber\\
M \cong M\tb \uob 
& = 
\bigoplus_x e_x M \otimes \C x ;\,\, m \mapsto \sum_x e_x m \otimes x.
\end{align}
The linear dual $M^*$ of $M$ has an $\FFF$-module structure,
which is determined by
$$
\langle e^x_y\, a\, m^\prime,\, m \rangle 
= 
\langle m^\prime,\, a^{-1}\, e^y_x m \rangle 
\quad (a \in G,\, x, y \in X,\,
m^\prime \in M^*,\, m \in M). 
$$ 
The module $M^*$ becomes a left dual object of $M$ via
\begin{align}
& 
ev \! :M^* \tb M \to \uob;\,\, 
\sum_x e_x m^\prime \otimes e^x m 
\mapsto 
\sum_y \langle e^y m^\prime, m\rangle y,
\nonumber\\
&
coev \! : \uob \to M \tb M^* ;\,\, 
x  
\mapsto 
\sum_i e^x m_i \otimes m^i,
\end{align}
where $\{ m_i \}$ denotes a basis of $M$ and 
$\{ m^i \}$ denotes its dual basis.
The canonical linear isomorphism $j_M \!:M \cong M^{**}$
becomes an isomorphism of $\FFF$-modules.
Hence ${\cal C}={}_{\FFF} {\bf mod}$ becomes a pivotal tensor category.

For each $M, N \in {\rm ob}\, {\cal C}$,
we define linear maps
$A_{M,N} \! : \CCC (\uob, M \tb N) \to \CCC (M^*, N)$,
$T_{M,N} \! : \CCC (M^*,N) \to \CCC (N^*,M)$, 
$E_{M,N} \! : \CCC (\uob, M \tb N) \to \CCC (\uob, N \tb M)$
by
\begin{gather*}
\begin{CD}
A_{M,N}(f)\!: M^* \cong M^* \tb \uob @> \id \tb f >> 
M^* \tb M \tb N 
@ > ev \tb \id>> \uob \tb N \cong N.
\end{CD}\\
T_{M,N} (g) = j_M^{-1} \circ g^*,
\quad
E_{M,N} (f) = (A_{N,M}^{-1} \circ T_{M,N} \circ A_{M,N}) (f),
\end{gather*}
respectively, where $f \in \CCC (\uob, M \tb N)$ and 
$g \in \CCC (M^*,N)$.

Then, the $r$-th indicator $\nu_r (M)$ of 
$M \in \mathrm{ob} {}_{\FFF} {\bf mod}$  
is defined by
$\nu_r (M):= \tr (E_{M,M^{\tb r-1}})$.

Let $M$ be a finite-dimensional vector space.
We say that $M$ is an {\it $\FFF$-space} if 
it is equipped with an associative action $\FFF \otimes M \to M$, 
that is,
the corresponding linear map $\pi_M\!: \FFF \to \mathrm{End} (M)$ satisfies
$\pi_M (\alpha\beta) = \pi_M (\alpha) \pi_M (\beta)$\,\,$(\alpha, \beta \in \FFF)$.   
Let $N$ be another $\FFF$-space. Then
$M \otimes N$ becomes an $\FFF$-space via
$\pi_{M\otimes N} (\alpha) = (\pi_M \otimes \pi_N)(\Delta (\alpha))$\,\,
$(\alpha \in \FFF)$.
For each $\FFF$-space $M$, we set
$\overline{M} := \pi_M (1) M$ and
$M^\FFF:= \pi_M (\inte) M$.
Then $\overline M$ becomes an $\FFF$-module.

\begin{lem}
\label{Fsp}
Let $M$ and $N$ be $\FFF$-spaces.\\ 
{\rm (1)}\, We have 
$\overline{\overline{M} \otimes N}
=
\overline{M \otimes N}
=
\overline{M \otimes \overline{N}}
$.\\
{\rm (2)}\,
Let $\cu^L$ be as in \eqref{eLdef}.
then, we have 
\begin{equation}
\label{MFFMLA}
M^\FFF = \{ m \in M\,|\, 
 \alpha m = \cu^L ( \alpha ) m
 \quad ( \alpha \in \FFF) \}.
\end{equation} 
{\rm (3)}\, 
The twist map $\mathrm{tw}_{M,N}\!: M\otimes N \to N\otimes M;$ 
$m\otimes n \mapsto n\otimes m$ 
satisfies 
$\mathrm{tw}_{M,N} \circ \pi_{M \otimes N} (\inte )= 
\pi_{N \otimes M} (\inte ) \circ \mathrm{tw}_{M,N}$.
In particular, it gives a linear isomorphism
$(M \otimes N)^\FFF \cong (N \otimes M)^\FFF$.
\end{lem}
\par\noindent
{\it Proof.}\,
Part (1) is obvious.
Let $N$ be the right-hand side of \eqref{MFFMLA}.
By \eqref{defint}, we have $M^\FFF \subseteq N$.
On the other hand, 
since $\cu^L (\inte) = 1$, we have
$n = \pi_M(\inte) n \in M^\FFF$ for each $n \in N$. 
Part (3) follows from
$$
(\mathrm{tw}_{M,N} \circ \pi_{M \otimes N} (\inte ))(m \otimes n)
= 
\frac{1}{|G|} \sum_{x,y,a} e^y_x a n \otimes e^x_y a m
=
(\pi_{N \otimes M} (\inte ) \circ \mathrm{tw}_{M,N})(m \otimes n).
$$
$\hfill \square$

\bigskip
For each $\FFF$-space $M$, 
there exists a linear isomorphism 
$\iota_{M} \!: M^\FFF \cong \CCC (\uob, \overline{M})$
such that
$\iota_{M} (m) (x) =e^x m$
for each $m \in M^\FFF$ and $x \in X$.
The inverse of $\iota_{M}$ is given by
$\iota_{M}^{-1} (f) =  \sum_{x \in X} f(x)$.

\begin{lem}
\label{CD}
For each $\FFF$-spaces $M$ and $N$, the diagram
$$
\begin{CD}
(M \otimes N)^\FFF @> \mathrm{tw}_{M,N} >> 
(N \otimes M)^\FFF\\
@ V \iota_{M \otimes N} VV @ VV \iota_{N \otimes M} V\\
\CCC (\uob, \overline{M \otimes N}) @> E_{\overline{M},\overline{N}} >> 
\CCC (\uob, \overline{N \otimes M})
\end{CD}
$$
is commutative.
\end{lem} 

\par\noindent
{\it Proof.}\,
Let $\sum_i m_i \otimes n_i$ be an element of $(M \otimes N)^\FFF$.
Set $g_1 = (A_{M,N} \circ \iota_{M \otimes N}) ( \sum_i m_i \otimes n_i )$
and $g_2 = (A_{N,M} \circ \iota_{N \otimes M}) ( \sum_i n_i \otimes m_i )$.
It is straightforward to verify that
$g_1 (m^\prime) = \sum_i \langle m^\prime,m_i \rangle
n_i$ for each $m^\prime \in {\overline{M}}^{\,*}$.
Hence 
$$
\langle m^\prime, T_{M,N} (g_1) (n^\prime )\rangle
= 
\sum_i \langle n^\prime, n_i \rangle \langle m^\prime, m_i \rangle
=
\langle m^\prime,  g_2 (n^\prime )\rangle
$$
for each $m^\prime \in \overline{M}^{\, *}$ and $n^\prime \in \overline{N}^{\,*}$.
This proves the assertion.
$\hfill \square$

\begin{prop}
For each
$M \in \mathrm{ob} {}_{\FFF} {\bf mod}$, we have
$\FS_r (M) = \nu_r (M)$. 
\end{prop}
\par\noindent
{\it Proof.}\,
Apllying Lemma \ref{CD} to $N = M^{\otimes r-1}$ and using Lemma \ref{Fsp} (1),
we obtain
$\nu_r (M) = \tr_{(M^{\otimes r})^\FFF} (\mathrm{tw}_{M,M^{\otimes r-1}})$.
Since $\inte$ is an idempotent, this equals to
$\tr \left(\pi_{M^{\otimes r}} (\inte ) \circ\mathrm{tw}_{M,M^{\otimes r-1}} \right) 
=
\tr \left(\pi_{M}^{\,\,\,\,\,\,\otimes r} (\Delta^{(r)} (\inte) ) \circ\mathrm{tw}_{M,M^{\otimes r-1}} \right) $ by Lemma \ref{Fsp} (3).
Hence, the assertion follows from the formula 
$
\tr ((f_1 \otimes \cdots \otimes f_r)\circ \mathrm{tw}_{M,M^{\otimes r-1}})
= \tr(f_1 \circ \cdots \circ f_r)  
$, which holds 
for each $f_1,\ldots, f_r \in \mathrm{End} (M)$.
$\hfill \square$

%\newpage

\end{document}